\documentclass[12pt]{amsart} %
\usepackage{amsfonts, amsmath, amsthm} %
\usepackage{graphicx} %
\usepackage{psfrag}

\newtheorem{pro}{Proposition}[section] %
\newtheorem{thm}[pro]{Theorem} %
\newtheorem{lem}[pro]{Lemma} %
\newtheorem{prop}[pro]{Property} %
\newtheorem{question}[pro]{Question} %
\newtheorem{example}[pro]{Example} %
\newtheorem{cnj}[pro]{Conjecture} %
\newtheorem{rmkk}[pro]{Remark} %
\newtheorem{assumptions}[pro]{Assumptions} %
\newtheorem{notation}[pro]{Notation} %

\newtheorem{clm}{Claim} %

\newtheorem{cor}[pro]{Corollary} %

\theoremstyle{definition} %
\newtheorem{dfns}[pro]{Definitions} %

\theoremstyle{remark} %

\newcommand{\pe}{Property $\mathcal{E}$} %
\newcommand{\pv}{Property $\mathcal{V}$} %
\newcommand{\pev}{Properties $\mathcal{V}$ and $\mathcal{E}$} %
\newcommand{\ccc}{counterclockwise cyclic order} %

\title{On Negami's Planar Cover Conjecture}

\date{\today} \address{Department of mathematical Sciences, University of
  Arkansas, Fayetteville, AR 72701}  
\address{Department of Mathematics, Nara Women's University Kitauoya
  Nishimachi, Nara 630-8506, Japan} 
\email{yoav@uark.edu} %
\email{yamasita@ics.nara-wu.ac.jp} %
\author{Yo'av Rieck} %
\author{Yasushi Yamashita} 

\def\marginpar#1{\mgp{\raggedright\tiny #1}}   \let\lbl=\label
\def\label#1{\lbl{#1}\ifinner\else\marginpar{\ref{#1} #1}\ignorespaces\fi}

\def \marginpar#1{}

\begin{document}


\date{\today}%
\begin{abstract}
Given a finite cover $f:\widetilde{G} \to G$ and an embedding of
$\widetilde{G}$ in the plane, Negami conjectures that $G$ embeds in $P^2$.
Negami proved this conjecture for regular covers.  In this paper we define two
properties (called \pev), depending on the cover $\widetilde{G}$ and its
embedding into $S^2$, and generalize Negami's result by showing: (1) If \pev\
are fulfilled then $G$ embeds in $P^2$.  (2) Regular covers always fulfill
\pev.  We give an example of an irregular cover fulfilling \pev.  Covers not
fulfilling \pev\ are discussed as well.
\end{abstract}

\maketitle

\section{Introduction}
\label{sec:intro}

In \cite{negami} S. Negami proved that if a graph $G$ has a finite,
unbranched, regular, planar cover then $G$ itself embeds in the projective
plane $P^2$ (for definitions see Section~\ref{sec:prelims}).  We call a graph
that embeds  in $P^2$ {\it projective}. Note that all planar graphs are
projective.  In the same paper, Negami conjectured that this holds in general:

\begin{cnj}[Negami's Conjecture]
\label{cnj:negami}
If a graph $G$ s has a finite unbranched planar cover then $G$ embeds in the
projective plane $P^2$.
\end{cnj}

\noindent Negami's result was extended
by S. Kitakubo to branched regular covers in \cite{kitakubo} (with the
exception of Section~\ref{sec:necessity}, any cover considered in this paper
may be branched).  We note that ``branched cover'' is not a standard term for
graphs; for a precise definition and discussion see
Definitions~\ref{dfn:branched-cover} and Remark~\ref{rmk:branched-cover}.  The
definition used by Kitakubo is what we call {\it weak cover}; it is immediate
from the definitions that every branched cover is a weak cover.  For regular
covers the converse holds as well: every regular weak cover is a branched
cover.  We do not have a clear idea about Negami's Conjecture for weak
covers:

\begin{question}
  Is Negami's Conjecture true for weak covers?
\end{question}

Given a finite cover $f:\widetilde{G} \to G$ and an embedding of $\widetilde{G}$ into
$S^2$ we define two properties called \pv\ and \pe.  These properties depend
on the covering map $f$ and the embedding of $\widetilde{G}$ into $S^2$.  We
prove Negami's Conjecture for covers $f:\widetilde{G} \to G$ fulfilling \pev.
We show that regular covers fulfill \pev\ (perhaps after re-embedding
$\widetilde{G}$ in $S^2$). The converse does not hold:  in
Example~\ref{ex:irregular} we show an irregular cover fulfilling \pev, showing
that the work here is more general than \cite{negami} and \cite{kitakubo}.  We
now give a more detailed description of our results and the structure of this
paper.

In Section~\ref{sec:prelims} we give the necessary definition and background.

In Section~\ref{sec:sufficiency} we define \pev\ (Definitions \ref{dfn:pv} and
\ref{dfn:pe}) and show (Theorem~\ref{thm:sufficiency} and
Corollary~\ref{cor:negami-cnj}) that if $f:\widetilde{G} \to G$ is a finite
planar cover fulfilling \pev\  then $G$ is projective, thus proving Negami's
Conjecture in that case. In fact, we show a little more: the map $f$ can be
extended to a map $f:S^2 \to F$, for some surface $F$, and this map is a
branched cover.  It is then easy to see (Lemma~\ref{pro:F-covered-bySorP})
that either $F \cong S^2$ or $F \cong P^2$. 

In Section~\ref{sec:regular} we show that if $\widetilde{G} \to G$ is a {\it
  regular} finite, planar cover then it fulfills \pev.

In \cite{negami} Negami mentions the following strategy for proving
Conjecture~\ref{cnj:negami}: given a finite cover $f:\widetilde{G} \to G$ it
is a well-known application of group theory that there exists finite cover
$\tilde{f} : \widetilde{\widetilde{G}} \to \widetilde{G}$, so that the
composition of the 
covers ({\it i.e.}, $f \circ \tilde{f} : \widetilde{\widetilde{G}} \to G$)  is
a finite regular cover (we remark that the degree of
$\widetilde{\widetilde{G}} \to \widetilde{G}$ is usually quite high). However,
Negami continues, even if we assume that $\widetilde{G}$ is planar, it does
{\it not} follow that $\widetilde{\widetilde{G}}$ is planar as well. In
Section~\ref{sec:necessity} we prove Theorem~\ref{thm:necessity}: let
$\widetilde{G} \to G$ be a finite {\it
unbranched} cover and let $\tilde{f}:\widetilde{\widetilde{G}} \to
\widetilde{G}$ be a cover with $\widetilde{\widetilde{G}} \subset S^2$.  If
the cover $f \circ \tilde{f}: \widetilde{\widetilde{G}} \to G$ fulfills \pev,
then either  $\widetilde{G}$ is planar and $f:\widetilde{G} \to G$ fulfills
\pev, or  $\widetilde{G}$ embeds in $P^2$; in that case, lifting
$\widetilde{G}$ to the universal cover of $P^2$ ({\it i.e.}, to the double
cover $S^2 \to P^2$) we obtain a finite planar cover of $G$ fulfilling \pev\
(this is sharp---see Example~\ref{ex:have-to-lift}).  Thus, if we
wish to pass to a cover of $\widetilde{G}$ in oder to prove that $G$ is
planar, we need not look for further than double covers.

In Section~\ref{sec:examples} we give examples.  The first is
Example~\ref{ex:irregular}, which is an irregular cover fulfilling \pev.  Next
(Example~\ref{ex:k4}) we give an example of two distinct double covers of the
planar graph $K_4$, both fulfilling \pev.  We use both covers to embed $K_4$.
The first gives an embedding into $S^2$ and the second into $P^2$,
demonstrating that the embedding depends on the cover in a non-trivial way.
Next, we exemplify Theorem~\ref{thm:necessity} by showing
(Example~\ref{ex:have-to-lift}) a planar cover $f:\widetilde{K}_4 \to K_4$
that does not fulfill \pev\ for any planar embedding of $\widetilde{K}_4$.
However, $\widetilde{K}_4$ embeds in $P^2$ and the lift of that embedding to
$S^2$ does fulfill \pev.  Finally, we show a more disturbing phenomenon: in
Conjecture~\ref{ex:no properties} we show a planar finite cover
$f:\widetilde{K}_4 \to K_4$ and we {\it conjecture} that is has no finite
cover $\tilde{f}:\widetilde{\widetilde{K}}_4 \to \widetilde{K}_4$ so the
composition $f \circ \tilde{f}:\widetilde{\widetilde{K}}_4 \to K_4$ fulfills
\pev. What happens here is that the cover $\widetilde{K}_4$ is simply the \em
wrong cover\em, and this cannot be fixed by passing to a higher cover. This,
perhaps, explains the difficulty in proving Negami's Conjecture.

So is Negami's Conjecture true?  At the current time, the answer is not known.
However, the work of many people (including D. Archdeacon, M. Fellows,
P. Hlin{\v{e}}n{\'y}, Negami, R. Thomas) accumulated in proving that Negami's
Conjecture is equivalent to the statement: {\it the graph $K_{1,2,2,2}$ has no
finite planar cover.}  This seems to us very strong evidence supporting
Conjecture~\ref{cnj:negami}.  \marginpar{CHECK} Assuming
Conjecture~\ref{cnj:negami} for a moment, in light of the results and examples
in this paper we ask: suppose we are given a finite planar cover
$f:\widetilde{G} \to G$ not fulfilling \pev\ for any embedding of
$\widetilde{G}$ (or a cover of $\widetilde{G}$) into $S^2$ ({\it i.e.}, a
``wrong cover''), how can we correct that cover?  We end
Section~\ref{sec:examples} with two ways to ``fix''  the cover in
Conjecture~\ref{ex:no  properties}.  The first is a cut-and-paste procedure
that produces a cover fulfilling \pev; this changes the graph $\widetilde{G}$.
The second does not change $\widetilde{G}$ but embeds it into some
non-orientable surface (not unlike Theorem~\ref{thm:necessity}).  This gives
an embedding of $G$ into {\it some} surface, not necessarily $P^2$, with
control over the Euler characteristic of that surface.

\section{Preliminaries}
\label{sec:prelims}

For a graph $G$, we denote the vertices of $G$ by $V(G)$ and the edges of $G$
by $E(G)$.  Naturally, any map between graphs is assumed to map vertices to
vertices and edges to edges.    We follow standard definitions and terminology
used in topology: $N(\cdot)$ means {\it closed} regular (or normal)
neighborhood, $\partial$ is read boundary, $\mbox{cl}$ is read closure, and
$\mbox{int}$ is read interior.  A homeomorphism is a continuous bijection with
a continuous inverse. $\chi(\cdot)$ stands for Euler characteristic.  All
surface are assumed to be connected.

As it is easy to reduce Negami's Conjecture to graphs $G$ with no cycles of
length one or two (that is, graphs $G$ in which every edge in $E(G)$ connects
distinct vertices and distinct edges have at most one vertex in common).  It
is also easy to reduce Negami's Conjecture to connected graphs $G$ and
$\widetilde{G}$.  Therefore, throughout this paper  we assume:

\begin{assumptions}
\label{assumptions} The graphs $G$ and $\widetilde{G}$ are connected and $G$
has no cycles of length 2 or less.
\end{assumptions}

We define branched covers in the two relevant situations, graphs and surfaces.
In these definitions disks are modeled on  $\{z \in \mathbb{C} : |z| < 1\}$,
and a {\it half disk} is a set homeomorphic to (and modeled on) $\{z \in
\mathbb{C} : |z| < 1 \mbox{ and } \Im(z) \geq 0 \}$ (here, $\Im(z)$ is the
imaginary part of $z$).  The boundary of a half disk are the points
corresponding to $\{z \in \mathbb{C} : \ |z| < 1 \mbox{ and } \Im(z) = 0\}$.
The symbol $\sqcup$ is used for disjoint unions.

\begin{dfns}
\label{dfn:branched-cover}
In~(1)--(3) below let $\widetilde{G}$ and $G$ be finite graphs and
$f:\widetilde{G} \to G$ a map.  In~(4)--(5) below let $F_1$ and $F_2$ be
compact surfaces and $f:F_1 \to F_2$ a proper map (that is, $f^{-1}(\partial
F_1) = \partial F_2$). 

\begin{enumerate}
\item   $f:\widetilde{G} \to G$ is called a {\it unbranched cover} if $f$ is
  onto and for any
  $\tilde{v} \in V(\widetilde{G})$ $f$ maps the neighbors of $\tilde{v}$
  bijectively onto the neighbors of $f(\tilde{v})$.
\item $f:\widetilde{G} \to G$  is called a {\it weak cover} if $f$ is onto and
  for any
  $\tilde{v} \in V(\widetilde{G})$ $f$ maps the neighbors of $\tilde{v}$ onto
  the neighbors of $f(\tilde{v})$.
\item $f:\widetilde{G} \to G$ is  called a {\it branched cover} if $f$ is onto
  and for any
  $\tilde{v} \in V(\widetilde{G})$ there is a positive integer $d =
  d(\tilde{v})$ so that every neighbor of $f(\tilde{v})$ has exactly $d$
  preimages that are neighbors of $\tilde{v}$, {\it i.e.}, the restriction of $f$
  to the neighbors of $\tilde{v}$ as a map to the neighbors of $f(\tilde{v})$
  is onto and $d$-to-1.  

  We call $d$ the {\it local
  degree} of $f$ at $\tilde{v}$.  If $d > 1$ then $\tilde{v}$ is called a {\it
  singular point}.  The set of all singular points is called the {\it singular
  set} and the image of the singular set is called the {\it branched set}.
\item map $f:F_1 \to F_2$ is called a {\it branched cover} if the
  following holds:
\begin{enumerate}
\item Every point $p \in \partial F_2$ has a neighborhood $D \ni p$ so that
  $D$ is half a  disk (with $p$ corresponding to 0) and $f^{-1}(D)$ is a
  disjoint collection of half disks $\sqcup_{i=1}^n D_i$ (for some $n$)  with
  $\partial D_i \subset \partial F_1$, so 
  that for every $i$, $f|_{D_i}:D_i \to D$ is a homeomorphism.
\item Every point $p \in \mbox{int}F_2$ has a neighborhood $D \ni p$ (with $p$
  corresponding to 0) so that $D$ is a disk and $f^{-1}(D)$ is a disjoint
  collection of disks $\sqcup_{i=1}^n D_i$ (for some $n$) so that for every
  $i$, $f|_{D_i}:D_i \to D$ is modeled on $z \mapsto z^d$ for some non-zero
  integer $d$.  

  If $f|_{D_i}$ is modeled on $z \mapsto z^d$ then $d$ is called the {\it
  local degree} at the center of $D_i$ (of course, for distinct values of $i$
  we may have distinct local degrees).  A point with local degree greater than
  one (in absolute value) is called a {\it singular point}, the union of the
  singular points is called the {\it singular set}, and the image of the
  singular set is called the {\it branched set}.  Note that the branch set is
  finite and contained in $\mbox{int}F_2$. 
\end{enumerate}
\end{enumerate}
\end{dfns}

Thus we see that every unbranched cover of graphs is a branched cover of
graphs with all local degrees one, and conversely a branched cover with all
local degrees one is an unbranched cover.  (Equivalently, unbranched covers
are covers with empty branch set). 

\begin{rmkk}
\label{rmk:branched-cover}{\rm
The definition of branched cover used in \cite{kitakubo} is the definition of
weak cover given above.  However, Kitakubo only considered regular covers.
It is left as an exercise to the reader to show that (under
Assumptions~\ref{assumptions}) a weak regular cover is in fact a branched cover.  We
do not know if Negami's Conjecture holds for weak covers, and it will probably
be a nice project to the reader to find and classify the counterexamples.
Weak covers will not appear in this paper again.}
\end{rmkk}

\begin{rmkk}{\rm
We identify $S^2$ with the Riemann sphere.  Then any rational function $f$
gives a branched cover $f:S^2 \to S^2$.  Conversely given a branched cover
$f:S^2 \to S^2$ we can multiply $f$ by $p^{-1}$ (for an appropriately chosen
polynomial $p$) so that no point in $\mathbb{C}$ is sent to $\infty$.  It then
follows from the Riemann Uniformization Theorem that (perhaps after
conjugation) $f/p$ is a polynomial.  Hence, after conjugating if necessary,
$f$ is a rational function.}
\end{rmkk}

We conclude this section with a few well-known lemmas about branched covers;
some of the proofs are sketched for the convenience of the reader.  For the
first lemma, the reader may consult, for example, \cite{hatcher} for the
definition of cover from the circle $S^1$ to itself.

\begin{lem}
\label{lem:cover-along-boundary}
Let $f:F_1 \to F_2$ be a branched cover between surfaces with non-empty
boundary.  Then the restriction $f|_{\partial F_1}$ an unbranched cover from
$f|_{\partial F_1}:\partial F_1 \to \partial F_2$.
\end{lem}

\begin{notation}
\label{notation:universal-cover}{\rm
The {\it universal cover } of $P^2$ is the map $\pi:S^2 \to P^2$ given by
identifying antipodal points.  It is an unbranched double cover.  }
\end{notation}

One of the basic facts about the universal cover of $P^2$ is:

\begin{lem}
\label{lem:covers-of-p2-factor}
Let $f:S^2 \to P^2$ be a cover.  Then $f$ factors through $\pi$, that is,
there exists a cover $f':S^2 \to S^2$ so that $\pi \circ f' = f$.
\end{lem}

\begin{proof}[Sketch of the proof] Let $B \subset P^2$ be the branch set.
Then $f|_{f^{-1}(P^2 \setminus B)}:f^{-1}(P^2 \setminus B) \to P^2 \setminus
B$ is an unbranched cover.  Since $P^2 \setminus B$ is non-orientable, it has
an {\it orientation double cover},\footnote{It is well-known that unbranched
covers correspond to subgroups of the fundamental group.  The orientation
double cover correspond to the group consisting of all the orientation preserving
loops in $P^2 \setminus B$.}  that is, an unbranched double cover $f_1:F \to
P^2 \setminus B$ from some orientable surface $F$.  A basic property of the
orientation double cover is that any cover from an orientable surface factors
through it.\footnote{A cover of $P^2 \setminus B$ is orientable if and only if
the corresponding subgroup contains only orientation preserving loops.  Any
such subgroup is contained in the subgroup of all orientation 
preserving loops, and hence any such cover factors through the orientation double
cover.}  Applying this in our setting, we see there a cover $f'|_{f^{-1}(P^2
\setminus B)}:f^{-1}(P^2 \setminus B) \to S^2$ so that $f'|_{f^{-1}(P^2
\setminus B)} \circ f_1 = f|_{f^{-1}(P^2 \setminus B)}$.  The surfaces $P^2
\setminus B$, $F$ and $f^{-1}(P^2 \setminus B)$ are not compact.  We
compactify them by adding one point to each end.  Denoting the
compactification of $F_1$ by $F_2$, we get a cover $f_2:F_2 \to P^2$ so that
$f$ factors through $f_2$.

All that remains to show is that $f_2$ is in fact the universal cover $\pi$.
We can easily see that $f_2$ is a double cover from a closed orientable
surface $F_2$ to $P^2$.  Since curves parallel to the punctures of $P^2
\setminus B$ are orientation preserving they lift to the double cover; hence
this cover is not branched.  Euler characteristic is multiplicative under
unbranched cover, and so $\chi(F_2) = 2$.  As $S^2$ is the only connected
surface with Euler characteristic 2, we get that $F_2 \cong S^2$.  By
uniqueness of the universal cover, $f_2 = \pi$ (perhaps after conjugation).
\end{proof}

The proof of the lemma below is an easy exercise in Euler characteristic:

\begin{lem}
\label{pro:F-covered-bySorP} If $S^2$ branch covers $F$ then
either $F \cong S^2$ or $F \cong P^2$.
\end{lem}

The following lemma tells us when $F$ is $P^2$; (3) is particularly convenient
since it requires only looking at one point:

\begin{lem}
\label{lem:when-is-F-on-orientable}

Let $F_1$, $F_2$ be surfaces and suppose $F_1$ is orientable.  Let $f:F_1 \to
F_2$ be a cover.  For $p \in F_2$, let $D_p$ be a open normal neighborhood of
$p$ endowed with and orientation (since $D_p$ is a disk this is always
possible).  For $q \in F_1$ with $f(q) = p$ let $D_q$ be a disk so that
$f|_{D_q}:D_q \to D_p$ is modeled on $z \to z^d$ (for some $d$).  Note
that the orientation on $D_p$ induces an orientation on $D_q$.

Then the following conditions are equivalent:

\begin{enumerate}
\item $F_2$ is non-orientable.
\item For any $p \in F_2$, there exist $q_1,q_2 \in F_1$ with $f(q_1) = f(q_2)
  = p$ so that the orientations induced orientations on $D_{q_1}$ and
  $D_{q_2}$ define opposite orientation on $F_1$.
\item For some $p \in F_2$ there exist $q_1,q_2 \in F_1$ with $f(q_1) = f(q_2)
  = p$ so that the orientations induced orientations on $D_{q_1}$ and
  $D_{q_2}$ define opposite orientation on $F_1$.
\end{enumerate}
\end{lem}

\begin{proof}
\noindent $(1) \Rightarrow (2)$: let $p \in
  F_2$ be an arbitrary point, and let $\gamma \subset F_2$ be an orientation
  reversing loop (which exists by assumption) based at $p$.  Since $F_1$
  admits no orientation reversing loops, the lift of $\gamma$ is an arc
  connecting two points (say $q_1$ and $q_2$) that project to $p$.  It is now
  easy to verify that the $D_{q_1}$ and $D_{q_2}$ induce opposite orientations
  on $F_1$ (in fact, this is equivalent to $\gamma$ being orientation
  reversing).  \medskip

\noindent $(2) \Rightarrow (3)$: Trivial.

\medskip

\noindent $(3) \Rightarrow (1)$: Given $p, q_1, q_2$ as in the
statement let $\alpha \subset F_1$ be any arc connecting $q_1$ and $q_2$.
Then the image of $\alpha$ on $F_2$ is an orientation reversing loop, and
hence $F_2$ is non-orientable.
\end{proof}

\section{\pv\ and \pe}
\label{sec:sufficiency}

Let $f:\widetilde{G} \to G$ be a finite planar cover. The embedding of
$\widetilde{G}$ into $S^2$ induces a cyclic order around each vertex of
$\widetilde{G}$. Our first condition, \pv, is a consistency condition
requiring that this order induces a
cyclic order around each vertex of $V(G)$.
Fix $v \in G(V)$ and $\tilde{v} \in v(\widetilde{G})$ in the preimage of  $v$.
There are two obstructions to inducing a cyclic order around $v$.  The first
obstruction is local: if $\tilde{v}$ is a singular vertex (say with local
degree $d$), then every neighbor of $v$ has $d$ preimages around $\tilde{v}$;
the ordering of the different preimages may contain inconsistencies.  For
example, if $v$ has three neighbors (say $v_1$, $v_2$, and $v_3$) and $d=2$
then there are six lifts of $v_1$, $v_2$ and $v_3$ adjacent to $\tilde{v}$,
say $\tilde{v}_1, \tilde{v}_2, \tilde{v}_3, \tilde{v}_1', \tilde{v}_2',
\tilde{v}_3'$ (with $f(\tilde{v}_i) = f(\tilde{v}_i') = v_i$, $i=1,2,3$).
These vertices may be cyclically ordered as $\tilde{v}_1, \tilde{v}_2,
\tilde{v}_3, \tilde{v}_1', \tilde{v}_3',  \tilde{v}_2'$. It is now not
possible to induce a cyclic order around $v$.  The second obstruction is
global: it is possible that the order around each preimage induces an order
around $v$ but different preimages induce distinct orders.  We now define \pv:

\begin{prop}[\pv]
\label{dfn:pv}  Let $f:\widetilde{G} \to G$ be a finite planar cover.  We say
that $f:\widetilde{G} \to G$ fulfills \pv\ if for any $v \in V(G)$ the following
two conditions are satisfied:

\begin{enumerate}
\item  For any $\tilde{v} \in f^{-1}(v)$, the cyclic order of the neighbors of
  $\tilde{v}$ (induced by the embedding into $S^2$) induces an order around
  $v$.  That is to say, denoting the neighbors of $v$ by $v_1,\dots,v_n$,
  after reordering the indices of if necessary, projecting the neighbors of
  $\tilde{v}$ to $V(G)$ in order we get $v_1, v_2, \dots, v_n,$  $v_1, v_2,
  \dots, v_n,\dots$ $v_1, v_2, \dots, v_n$.   (This condition is vacuous is
  $\tilde{v}$ is not singular).
\item The order obtained is independent of choice of preimage.
\end{enumerate}
\end{prop}

\begin{figure}
\psfrag{v1}{$v_1$}
\psfrag{v2}{$v_2$}
\psfrag{v3}{$v_3$}
\psfrag{v4}{$v_4$}
\psfrag{vt}{$v_t$}
\psfrag{vpt}{$v_t'$}
\centerline{\includegraphics[height=3cm]{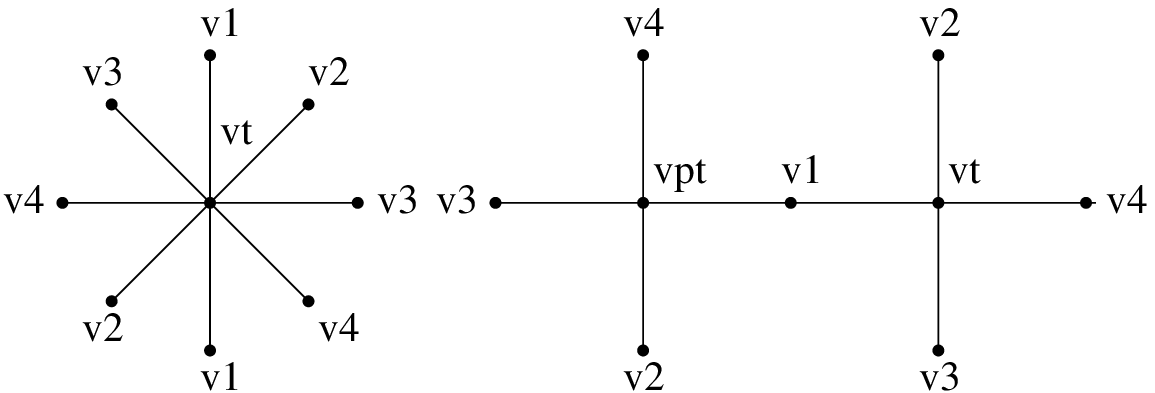} }
\caption{Covers not fulfilling \pv}
\label{fig:no-pv}
\end{figure}

\begin{figure}
\psfrag{v1}{$v_1$}
\psfrag{v2}{$v_2$}
\psfrag{v3}{$v_3$}
\psfrag{v4}{$v_4$}
\psfrag{vt}{$v_t$}
\psfrag{vpt}{$v_t'$}
\centerline{\includegraphics[height=3cm]{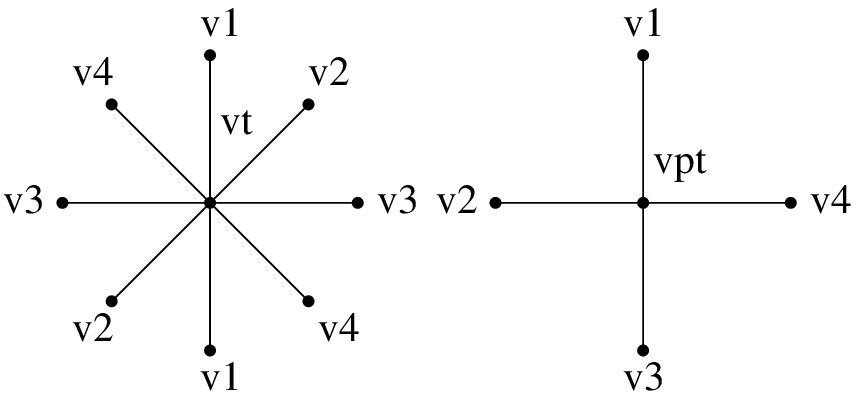} }
\caption{A cover fulfilling \pv}
\label{fig:yes-pv}
\end{figure}

\noindent {\it Example:} Suppose that $v \in V(G)$ has four
neighbors. Figure~\ref{fig:no-pv} shows the two possibilities that a cover
does not fulfill \pv.  The vertices $\tilde{v},\tilde{v}'$ in that figure
project to $v$; their neighbors are labeled using the labels
$v_1.v_2,v_3,v_4$, according to which vertex of $G$ they project to.  Using
the same labeling scheme, in Figure~\ref{fig:yes-pv} we show a little part of
a cover that fulfills \pv.


Let $f:\widetilde{G} \to G$ be a finite planar cover fulfilling \pv. Given $v
\in V(G)$ and $\tilde{v} \in f^{-1}(v)$, moving around a the neighbors of
$\tilde{v}$ counterclockwise allows us to distinguish a particular cyclic
order from its reverse order
on the neighbors of $v$.  We call this the {\it \ccc. }  
We assign $\tilde{v}$ an arbitrary sign (plus or minus).\footnote{The freedom
  of choice will be useful in Section~\ref{sec:necessity}.}  
If $\tilde{v}'$ is another preimage of $v$
we assign  $\tilde{v}'$ the same sign as $\tilde{v}$ if the \ccc\ around
$\tilde{v}'$ projects to the same order around $v$, and the opposite sign
otherwise. A sign assignment as above is called {\it valid}. We refer the
reader to Section~\ref{sec:examples} for examples of valid sign assignments. 
Like \pv, \pe\ is a question of consistency:

\begin{prop}[\pe]
\label{dfn:pe} Let $f:\widetilde{G} \to G$ be a finite planar
cover fulfilling \pv\ with valid signs on $G(\widetilde{V})$.  We say that
$f:\widetilde{G} \to G$ fulfills \pe\ if for any edge $e \in E(G)$ and any two
edges $\tilde{e}_1, \tilde{e}_2 \in E(\widetilde{G})$ that project to $e$, we
have that $\tilde{e}_1$ connects vertices of the same sign if and only if
$\tilde{e}_2$ does.
\end{prop}

We now state and prove Theorem~\ref{thm:sufficiency}. In that theorem, we
consider a graph $\widetilde{G}$ embedded in $S^2$ ($G$ embedded in some
surface $F$, resp.).  The closure of the components of $S^2 \setminus
\widetilde{G}$ ($F \setminus G$, resp.) are called {\it faces} of $S^2$ ($F$,
resp.). 

\begin{thm}
\label{thm:sufficiency} 
Let $\widetilde{G}$ be a connected planar graph, and let $f:\widetilde{G} \to
G$ be a finite cover.  Then $f:\widetilde{G} \to G$ fulfills \pev\ if and only
if there exist a surface $F$ containing $G$ and a map $f':S^2 \to F$ with the
following properties: 

\begin{enumerate}
\item $f'$ extends $f$, that is, for every point $p \in \widetilde{G}$, $f'(p)
  =f(p)$.

\item $f'$ is a branched cover, and hence (by
  Lemma~\ref{pro:F-covered-bySorP}) $f\cong P^2$ or $F \cong S^2$.

\item The intersection of the branch points of $f'$ with $G$ is contained in
  $V(G)$.  More specifically, it is exactly the branch set of $f$.

\item The intersection of the singular points of $f'$ with $\widetilde{G}$ is
  contained in $V(\widetilde{G})$.  More specifically, it is exactly the set
  of singular points of $f$ and with the same local degrees.

\item The faces of $S^2$ ($F$, resp.) are all disks, and every such face
  contains at most one singular point (branch point, resp.).

\item\label{order} For each $v \in V(G)$, the cyclic order induced on the
  neighbors of $v$\ by $f$ and the cyclic order given by the embedding of $G$
  into $F$ are the
  same.
\end{enumerate}
\end{thm}

Negami's Conjecture for finite planar covers fulfilling \pev\
(Corollary~\ref{cor:negami-cnj} below) follows easily from points (1) and (2)
above.  However, since Negami's Conjecture does not require (3)--(6),
corollary~\ref{cor:negami-cnj} is {\it not} and ``if and only if'' statement.

\begin{cor}[Negami's Conjecture for covers fulfilling \pev]
\label{cor:negami-cnj}%
Let $\widetilde{G}$ be a connected graph, $\widetilde{G} \subset S^2$. Let
$f:\widetilde{G} \to G$ be a cover fulfilling \pev. Then $G$ is projective.
\end{cor}

\begin{proof}[Proof of Theorem~\ref{thm:sufficiency}]
Suppose $f:\widetilde{G} \to G$ is a finite planar cover fulfilling \pev. We
construct the surface $F$ and the map $f':S^2 \to F$:

\smallskip

\noindent {\bf Step One: Vertices.}  We first construct $F$ near each vertex
of $G$ and extend the map $f'$ from a neighborhood of $V(\widetilde{G})$ to
this neighborhood of $V(G)$.  For each $\tilde{v} \in V(\widetilde{G})$ let
$D_{\tilde{v}}$ be a regular neighborhood of $\tilde{v}$.  By the Normal
Neighborhood Theorem we may assume each $D_{\tilde{v}}$ is a closed ({\it
sic.}) disk and these disks are disjoint; moreover, the edges of
$\widetilde{G}$ intersect $D_{\tilde{v}}$ in radial arcs, as in
Figure~\ref{fig:D-v}.
\begin{figure}
\psfrag{Dvt}{}
\centerline{\includegraphics[height=3cm]{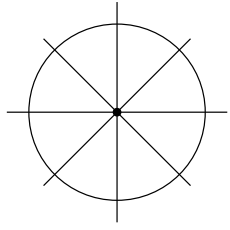} }
\label{fig:D-v}
\caption{$D_{\tilde{v}}$}
\end{figure}
For $v \in V(G)$ we pick a small disk $D_v$ and embed a neighborhood of $v$ in
$D_v$ as follows: $v$ is the center of $D_v$, and $E(G)$ intersects $D_v$ in
radial arcs, each arc corresponding to the tip of an edge that is incident to
$v$.  This can be done in several different ways resulting in distinct cyclic
orders around $v$.  We embed $v$ in $D_v$ so that the cyclic order around $v$
agrees with the order induced by the neighbors of the preimages of $v$; this
is well-defined since the cover fulfills \pv.

We assign $V(\widetilde{G})$ a valid sign convention.
Endowing $D_v$ with an orientation gives a specific \ccc\ around $v$. We
orient $D_v$ so that the counterclockwise cyclic order around $D_v$ agrees
with the \ccc\ along induced by the neighbors of a vertex of $f^{-1}(v)$ with
plus sign and is opposite the \ccc\ induced by vertices of negative sign.
Since the sign assignment is valid, this is well defined.

The surface constructed so far, $\sqcup_{v \in V(G)} D_v$, is denoted
$F_V$. Next we extend the map $f$ to $f':\sqcup_{\tilde{v} \in
V(\widetilde{G})} \to F_V$ by mapping $D_{\tilde{v}}$ to $D_{f(\tilde{v})}$ by
a map modeled $z \mapsto z^d$ (where $d$ is the local degree at $\tilde{v}$).

We check each property listed in Theorem~\ref{thm:sufficiency}, in the same
order:

\begin{enumerate}
\item  Holds by construction: for any point $p \in
  \widetilde{G} \cap (\sqcup_{\tilde{v} \in V(\widetilde{G})} D_{\tilde{v}})$,
  $f'(p) = f(p)$.
\item Again by construction, $f':\sqcup_{\tilde{v} \in V(\widetilde{G})}
  D_{\tilde{v}} \to F_V$ is a branched cover.
\item  Holds by construction.  
\item  Holds by construction. 
\item  To be verified later.
\item  Holds by construction.  Since this is a local property and we will not
  modify $f'$ near $V(\widetilde{G})$ any more, we will not need to check this
  property again.
\end{enumerate}

\begin{rmkk}
\label{rmk:orientation}
{\rm
$f':D_{\tilde{v}} \to D_v$ is orientation preserving if and only if the sign
of $v$ is plus.}
\end{rmkk}

\smallskip

\noindent {\bf Step Two: Edges.} Next, we extend the construction of $F_V$
to a neighborhood of $E(G)$ and extend the range of $f'$ from a neighborhood
of $E(\widetilde{G})$.    The neighborhood of an edge $e \in E(G)$ (say $e =
(v,u)$, for some $v,u \in V(G)$) is a closed ({\it sic.}) band that connects
$D_v$ to $D_u$.  We first glue the band to $D_v$ and extend the orientation of
$D_v$ along the band.  We then glue the band to $D_u$ as follows: choose an
edge $\tilde{e} \in f^{-1}(e)$.  If $\tilde{e}$ 
connects vertices with the same sign, the orientation of the band agrees
with the orientation of $D_u$ and if $\tilde{e}$ connects vertices with
the opposite signs, the orientation of the band disagrees with the orientation
of $D_u$.\footnote{Since the disks $D_v$ are oriented, we can place them on a
coffee table all facing up.  For the edge $e = (u,v)$, the band $N(e)$
connects $D_u$ to $D_v$.  If $\tilde{e}$ connects vertices of the same
sign, this band is untwisted ({\i.e.}, lies flat on the table) and if
$\tilde{e}$ connects vertices of the opposite sign, the band is twisted.}
Since $f$ fulfills \pe\ this construction is independent of choice.  We denote
the part of $F$ constructed so far $F_E$.

Again, we check each property listed in Theorem~\ref{thm:sufficiency},
ignoring the properties we are done with.

\begin{enumerate}
\item  Holds by construction for any point of $\widetilde{G}$. Since this is a
local property and we will not modify $f'$ on $(f')^{-1}(F_E)$ any more, we
will not need to check this property again.
\item  The map $f':(f')^{-1}(F_E) \to F_E$ is a branched cover by construction.
\item and (4) Since we did not introduce any new branch points or singular points,
(3) and (4) still hold.  We will not change $f'$ near $\widetilde{G}$ so do not need
to check these properties again,
\setcounter{enumi}{4}
\item  Note that so far, all singular (resp.  branch) points are vertices of
$\widetilde{G}$ (resp. $G$). (5) will be verified later.
\end{enumerate}

\smallskip

\noindent {\bf Step Three: Closing $F_E$.}  We have constructed a
branched cover  $f':(f')^{-1} (F_E) \to F_E$, both compact surfaces with
non-empty boundary.  Let $\tilde{\gamma}$ denote a boundary component of
$f^{-1} (F_E)$.  By Lemma~\ref{lem:cover-along-boundary} $f|_{f^{-1}(F_E)}$
maps $\tilde{\gamma}$ to a boundary component of $F_E$, say $\gamma$, and the
map $f|_{\tilde{\gamma}}:\tilde{\gamma} \to \gamma$ is a cover.  Since
$\tilde{\gamma}$ and $\gamma$ are both circle, $f|_{\tilde{\gamma}}$ is
modeled on the restriction of $z \mapsto z^d$ to the unit circle  (for some
$0 \neq d \in \mathbb Z$, called the {\it winding number } of
$f|_{\tilde{\gamma}}$). Since $\widetilde{G}$ is connected and $f^{-1}(F_E)$
is a neighborhood of $\widetilde{G}$, ${f^{-1}(F_E)}$ is connected as well;
hence the components of $S^2 \setminus f^{-1}(F_E)$ are all disks.  Let
$D_{\tilde{\gamma}}$ be the closed disk bound by $\tilde{\gamma}$ disjoint
from $\widetilde{G}$, that is, the closure of the component of $S^2 \setminus
{f^{-1}(F_E)}$ adjacent to $\tilde{\gamma}$.  We attach a disk (say
$D_{\gamma}$) to $\gamma$ and extend the map $f'$ by mapping
$D_{\tilde{\gamma}}$ to $D_{\gamma}$ by coning.\footnote{Recall that we model
both $D_{\tilde{\gamma}}$ and $D_\gamma$ on the unit disk in $\mathbb C$.  By
{\it coning} we mean that $f'$ is modeled on $z \to z^d$, extending
$f'_{\tilde{\gamma}}$.} If the absolute value of the winding number
$f|_{\tilde{\gamma}}$ is more than one we introduce exactly one singular point
on $D_{\tilde{\gamma}}$ and one branch point on $D_{\gamma}$, otherwise no new
singular or branch point is introduced. Continuing this way, we cap off every
component of the boundary of $F_E$, finally constructing a closed surface $F$
containing $G$ and a branched cover $f:S^2 \to F$.  With the exception of  (2)
and (5), we verified that the cover $f':S^2 \to F$ fulfills all the conditions
of Theorem~\ref{thm:sufficiency}.  We note that $f'|_{f^{-1}(F_E)}$ is a
branched cover, and by construction, $f'|_{\mbox{cl}(F \setminus F_E)}$ is a
branched cover as well; this establishes (2).  It is straightforward to see
that the construction gives (5) as well.

Conversely, given $f:\widetilde{G} \to G$ and a map $f'$ as in the statement
of Theorem~\ref{thm:sufficiency}, the embedding of $G$ into $F$ induces a
cyclic order around the vertices of $G$.  Lifting these orders to
$\widetilde{G}$ we get the order around each vertex of $\widetilde{G}$; this
order coincides with the order given by the embedding of $\widetilde{G}$ into
$S^2$.  It is now easy to see that \pev\ follow.

This completes the proof of Theorem~\ref{thm:sufficiency}.
\end{proof}

We can know if $F \cong S^2$ or $F \cong P^2$ by looking at labels only:

\begin{pro}
\label{pro:p2iff-both-signs}
The following are equivalent:
\begin{enumerate}
\item The surface $F$ constructed in Theorem~\ref{thm:sufficiency} is
homeomorphic to $P^2$.
\item There exists a vertex $v \in V(G)$ with preimages of opposite signs.
\item Every vertex $v \in V(G)$ has preimages of opposite signs.
\end{enumerate}
\end{pro}

\begin{proof}
This is immediate from Remark~\ref{rmk:orientation} and
Lemma~\ref{lem:when-is-F-on-orientable}.
\end{proof}

\section{Regular covers}
\label{sec:regular}

In this section we show that if $\widetilde{G} \to G$ is a finite planar
regular cover then it fulfills \pev.  Some of the material in this section is
from \cite{negami} and \cite{kitakubo}.

\begin{pro}
\label{pro:regular-is-nice} Let $f:\widetilde{G} \to G$ be a
regular planar cover.  Then (perhaps after re-embedding $\widetilde{G}$)
$f:\widetilde{G} \to G$ fulfills \pev.
\end{pro}

\begin{proof}

Recall that a cover $f:\widetilde{G} \to G$ is called {\it regular} if there
is a group $\Gamma$ acting  on $\widetilde{G}$ so that for any
$\tilde{v},\, \tilde{v}' \in V(\widetilde{G})$ with $f(\tilde{v}) =
f(\tilde{v}')$ there exist $\gamma \in \Gamma$ with $\gamma(\tilde{v}) =
\tilde{v}'$, and for any $\tilde{e},\, \tilde{e}' \in E(\widetilde{G})$ with
$f(\tilde{e}) = f(\tilde{e}')$ there exist $\gamma \in \Gamma$ with
$\gamma(\tilde{e}) = \tilde{e}'$.

\begin{lem}
\label{lem:translations-extend} Let $f:\widetilde{G} \to G$ be a
finite planar regular cover with group $\Gamma$.  Then (perhaps after
re-embedding $\widetilde{G}$) the action of $\Gamma$ can be extended to an
action on $S^2$.
\end{lem}

\begin{proof}[Proof of Lemma~\ref{lem:translations-extend}]
The proof is an induction on the number of vertices of $\widetilde{G}$; when
$\widetilde{G}$ is not 3-connected we reduce this number.  Hence the base case
of the induction is:

\medskip

\noindent{\bf Base Case:} $\widetilde{G}$ is 3-connected.  This is
a well-known theorem of Whitney \cite{whitney}.

\medskip

\noindent{\bf Inductive Step:} We assume $\widetilde{G}$ is not
3-connected. The proofs for 2-connected and 1-connected graphs are similar,
and we omit the easier case of 1-connected graphs.  Let $\{ \tilde{c}_1,
\tilde{c}_2 \}$ be a cut pair for $\widetilde{G}$ and suppose removing
$\tilde{c}_1,\tilde{c}_2$ from $\widetilde{G}$ we obtain the graphs
$\widetilde{G}_0$ and $\widetilde{G}_1$; we assume further that $\{
\tilde{c}_1,\tilde{c}_2 \}$ were chosen so that $\widetilde{G}_0$ is minimal
with respect to inclusion.\footnote{This implies that $\widetilde{G}_0$ is
2-connected, but $\widetilde{G}_1$ need not be connected.}  Note that
$V(\widetilde{G}_0) \neq \emptyset$ and $V(\widetilde{G}_1) \neq \emptyset$.
For any $\gamma \in \Gamma$, $\{\gamma(\tilde{c}_1),\gamma(\tilde{c}_2) \}$ is
a cut pair.  Therefore minimality of  $\widetilde{G}_0$ implies that
$\gamma(\widetilde{c}_i) \not\in  \widetilde{G}_0$ ($i=1,2$).  We use the
notation $\mbox{ext}(\widetilde{G}_0)$ (called the {\it extension} of
$\widetilde{G}_0$) for the subgraph spun by $V(\widetilde{G}_0) \cup
\{\tilde{c}_1,\tilde{c}_2 \}$ and similarly $\mbox{ext}(\widetilde{G}_1)$ is
the graph spun by $V(\widetilde{G}_1) \cup \{\tilde{c}_1,\tilde{c}_2
\}$. Since $\widetilde{G}$ is 2-connected there is a pass (say $\alpha_1$) in
$\widetilde{G} \setminus \widetilde{G}_0$ (equivalently, in
$\mbox{ext}(\widetilde{G}_1)$) connecting $\tilde{c}_1$ to $\tilde{c}_2$.  By
choosing the shortest pass, we can guarantee that $\alpha_1$ is embedded.
Therefore, the graph obtained by adding an edge (say $\tilde{e}$) that
connects $\tilde{c}_1$ to $\tilde{c}_2$ to $\mbox{ext}(\widetilde{G}_0)$ is
planar.   We denote this graph by $\mbox{ext}(\widetilde{G}_0) \cup \{
\tilde{e} \}$.

Let $\Gamma_0$ be the (possibly trivial) subgroup of $\Gamma$ that leaves
$\mbox{ext}(\widetilde{G}_0)$ invariant (equivalently, leaves
$\widetilde{G}_0$ invariant).  Then for any $\gamma \in \Gamma_0$,
$\{\gamma(\tilde{c}_1),\gamma(\tilde{c}_2) \} \subset
\mbox{ext}(\widetilde{G}_0)$, and hence
$\{\gamma(\tilde{c}_1),\gamma(\tilde{c}_2) \} = \{\tilde{c}_1,\tilde{c}_2 \}$.
Therefore the action of $\Gamma_0$ can be extended to an action on
$\mbox{ext}(\widetilde{G}_0) \cup \{ \tilde{e} \}$.  Since $V(\widetilde{G}_1)
\neq \emptyset$, $|V(\mbox{ext}(\widetilde{G}_0) \cup \{ \tilde{e} \})| <
|V(\widetilde{G})|$ (where $|V(\cdot)|$ denotes number of vertices) and we may
apply the inductive hypothesis to get a re-embedding of
$\mbox{ext}(\widetilde{G}_0) \cup \{ \tilde{e} \}$ into $S^2$ so that the
action of $\Gamma_0$ extends to $S^2$.  Removing $\tilde{e}$, we obtain a
re-embedding of $\mbox{ext}(\widetilde{G}_0)$, denoted $\widehat{G}_0$.

\begin{rmkk}
\label{rmk:hat-is-same} {\rm%
It is not hard to argue that $\mbox{ext}(\widetilde{G}_0) \cup \{ e \}$ is
3-connected.  Hence by Whitney \cite{whitney} (the base case of the induction)
$\widehat{G}_0$ is the original embedding of $\mbox{ext}(\widetilde{G}_0)$.}
\end{rmkk}

Since $\widetilde{G}$ is 2-connected there is a pass (say $\tilde{\alpha}_0$)
in $\mbox{ext}(\widetilde{G}_0)$ connecting $\tilde{c}_1$ to $\tilde{c}_2$.
By choosing the shortest pass, we can guarantee that $\tilde{\alpha}_0$ is
embedded.  Similar to the process above, we replace $\widetilde{G}_0$ with an
edge (say $\tilde{e}_0$) connecting $\tilde{c}_1$ and $\tilde{c}_2$, obtaining
the planar graph $(\widetilde{G} \setminus \widetilde{G}_0) \cup \tilde{e}$.
Of course, we cannot expect $\Gamma$ to act on this graph.  To that end, we
repeat this operation on the image of $\widetilde{G}_0$ under $\Gamma$ and
obtain the planar graph:
$$(\widetilde{G} \setminus \cup_{\gamma \in \Gamma}
\gamma(\widetilde{G}_0)) \cup (\cup_{\gamma \in \Gamma} \gamma(\tilde{e}_0)).$$%
Denote this graph by $\widetilde{H}$, and note that by construction $\Gamma$
acts on $\widetilde{H}$. Since $V(\widetilde{G}_0) \neq
\emptyset$,  $|V(\widetilde{H})| < |V(\widetilde{G})|$ and we may apply the
inductive hypothesis to get a re-embedding of $\widetilde{H}$  into $S^2$ so
that the action of $\Gamma$ extends to $S^2$.  For every $\gamma \in \Gamma$,
We replace $\tilde{e}_0$ by $\gamma(\widehat{G}_0)$.
It is now easy to see that we obtain an embedding of $\widetilde{G}$ into
$S^2$ and the action of $\Gamma$ on this graph extends to $S^2$, as desired.

This completes the proof of Lemma~\ref{lem:translations-extend}.
\end{proof}

After extending the action of $\Gamma$ to $S^2$, we denote the group elements
by $\phi:S^2 \to S^2$.

\begin{lem}
\label{lem:extends-implies-pev}
Let $f:\widetilde{G} \to G$ be a finite planar regular cover with group
$\Gamma$.  If any $\gamma \in \Gamma$ can be extended to a homeomorphism
$\phi:S^2 \to S^2$ then the cover fulfills \pev.
\end{lem}

\begin{proof}[Proof of Lemma~\ref{lem:extends-implies-pev}]

First,  let $\tilde{v}$ be a singular vertex (say $\tilde{v}$ projects to
$v$). We need to show that the cyclic order of the neighbors of $\tilde{v}$
induces an order on the neighbors $v$.

\begin{clm}
Let $\tilde{v}_1, \tilde{v}_2 \in V(\widetilde{G})$ be consecutive vertices in
the \ccc\ around $\tilde{v}$ and let $\tilde{v}_1'$ be a neighbor of
$\tilde{v}$ so that $f(\tilde{v}_1') = f(\tilde{v}_1)$

Then the vertex that follows $\tilde{v}_1'$ in the \ccc\ around $\tilde{v}$
projects to the same vertex as $\tilde{v}_2$.
\end{clm}

\noindent Proof of Claim~1: Since $\tilde{v}_1$ and $\tilde{v}_1'$
project to the same vertex, there exists $\gamma \in \Gamma$ so that
$\gamma(\tilde{v_1}) = \tilde{v}_1'$. Denote the edge
$\{\tilde{v},\tilde{v}_1'\}$ by $e_1$ and the image of
$\{\tilde{v}_1,\tilde{v}\}$ under $\gamma$ by $e_2$ (that is,
$e_2=\{\tilde{v_1}',\gamma(\tilde{v})\}$).  If $\tilde{v} \neq
\gamma(\tilde{v})$ then $e_1 \cup e_2$ projects to a cycle of length 1 or 2,
contradicting our assumption (recall Assumptions~\ref{assumptions}).  We
conclude that $\gamma$ fixes $\tilde{v}$. By assumption, there exists
$\phi':S^2 \to S^2$ a homeomorphism that extends $\gamma$;  thus
$\phi(\tilde{v}) = \tilde{v}$, $\phi(\tilde{v_1}) = \tilde{v_1}'$, and
$\phi(\tilde{v_2})$ is a vertex that projects to the same vertex as
$\tilde{v}_2$, and follows $\tilde{v}_1'$ in the \ccc\ around
$\tilde{v}$. This proves Claim~1.

It follows immediately from Claim~1 that the \ccc\ around $\tilde{v}$ indices
a cyclic order around $v$.

Next, let $\tilde{v}, \tilde{v}' \in V(\widetilde{G})$ be distinct vertices
that project to the same vertex $v \in V(G)$. Then there exists $\gamma \in
\Gamma$ so that $\gamma(\tilde{v}) = \tilde{v}'$, and by assumption there
exists $\phi:S^2 \to S^2$ extending $\gamma$.  It is easy to see that $\phi$
induces an order preserving bijection between the neighbors of $\tilde{v}$ and
the neighbors of $\tilde{v}'$.  Therefore the neighbors of $\tilde{v}$ and the
neighbors of $\tilde{v}'$ induce the same order on the neighbors of $v$; this
establishes \pv.

Since $f:\widetilde{G} \to G$ fulfills \pv, we can assign a valid sign assignment
for $V(\widetilde{G})$ (as described in Section~\ref{sec:sufficiency}). Let
$\tilde{e}_1, \tilde{e}_2 \in E(\widetilde{G})$ be edges that project to the
same edge $e \in E(G)$; say $\tilde{e}_1 = \{\tilde{v}_1,\tilde{v}_1'\}$ and
$\tilde{e}_2 = \{\tilde{v}_2,\tilde{v}_2'\}$. Than there exists $\gamma \in
\Gamma$ so that $\gamma(\tilde{e}_1) = \tilde{e}_2$, equivalently
$\{\gamma(\tilde{v}_1),\gamma(\tilde{v}_2)\} = \{\tilde{v}_1',\tilde{v}_2'
\}$.  If $\gamma$ is preserves the orientation of $S^2$, the sign of
$\tilde{v}_1$ is the same as the sign of $\gamma(\tilde{v}_1)$ and the sign of
$\tilde{v}_2$ is the same as the sign of $\gamma(\tilde{v}_2)$; if $\gamma$
reverses the orientation of $S^2$ the sign of $\tilde{v}_1$ is the opposite
the sign of $\gamma(\tilde{v}_1)$ and the sign of $\tilde{v}_2$ is opposite
the sign of $\gamma(\tilde{v}_2)$.  In both cases, $\tilde{e}_1$ connects
vertices of the same sign if and only if $\tilde{e}_2$ does. This establishes
\pe\  and completes the proof of Lemma~\ref{lem:extends-implies-pev}.
\end{proof}

Clearly, Proposition~\ref{pro:regular-is-nice} follows from
Lemmas~\ref{lem:translations-extend} and \ref{lem:extends-implies-pev}.
\end{proof}

\section{Higher covers}
\label{sec:necessity}

In this section, we try to better our situation by passing to higher covers.
Let $f:\widetilde{G} \to G$ be a finite cover.  When does there exist a finite
planar cover $\tilde{f}:\widetilde{\widetilde{G}} \to \widetilde{G}$ so that
the composition $\tilde{f} \circ f: \widetilde{\widetilde{G}} \to G$ fulfills
\pev?  Assume such cover exists.  To address this question, we wish to apply
Theorem~\ref{thm:sufficiency} to the cover
$\tilde{f}:\widetilde{\widetilde{G}} \to \widetilde{G}$.  However, it does not
follow that $\tilde{f}:\widetilde{\widetilde{G}} \to \widetilde{G}$ fulfills
\pev.  To guarantee that, we need to add the assumption that the cover
$f:\widetilde{G} \to G$ is not branched.

Recall from Notation~\ref{notation:universal-cover} that the map $\pi:S^2
\to P^2$ given by identifying antipodal points is an unbranched, double cover
called the {\it universal cover} of $P^2$.  Therefore, given a graph
$\widetilde{G} \subset P^2$, $\pi^{-1}(\widetilde{G})$ is a graph that double
covers $\widetilde{G}$ and naturally embeds in $S^2$  

\begin{thm}
\label{thm:necessity} Let $f:\widetilde{G} \to G$ be an unbranched
finite (not necessarily planar) cover and $\tilde{f}:\widetilde{\widetilde{G}}
\to \widetilde{G}$  be a finite planar covers. Suppose that the composition
$\tilde{f} \circ f:\widetilde{\widetilde{G}} \to G$ fulfills \pev.  Then one
of the following holds:

\begin{enumerate}
\item  $\widetilde{G}$ is planar and for some embedding of $\widetilde{G}$
  into $S^2$, $f:\widetilde{G} \to G$ fulfills \pev.
\item  $\widetilde{G}$ embeds in $P^2$ so that the cover $f \circ \pi:
  \pi^{-1}(\widetilde{G}) \to G$ fulfills \pev.
\end{enumerate}
\end{thm}

\begin{proof}
We begin by showing:

\begin{pro}
\label{pro:tilde-f-has-pev}
Let $f:\widetilde{G} \to G$ be an unbranched finite planar cover and
$\tilde{f}:\widetilde{\widetilde{G}} \to \widetilde{G}$  be a finite planar
covers.  Suppose that the composition $\tilde{f} \circ
f:\widetilde{\widetilde{G}} \to G$ fulfills \pev.  Then
$\tilde{f}:\widetilde{\widetilde{G}} \to \widetilde{G}$ fulfills \pev.
\end{pro}

\begin{proof}[Proof of Proposition~\ref{pro:tilde-f-has-pev}]
Let $\tilde{\tilde{v}}, \tilde{\tilde{v}}' \in V({\widetilde{\widetilde{G}}})$
be two vertices that project to the same vertex under $\tilde{f}$ (say
$\tilde{v}$) and denote $f(\tilde{v})$ by $v$. By assumption,  $\tilde{f}
\circ f$ fulfills \pv.  Therefore, the cyclic order on the neighbors
$(\tilde{f} \circ f)^{-1}(v)$ induces a cyclic order on the neighbors of
$v$. Since $f$ is unbranched, the restriction of $f$ to the neighbors of
$\tilde{v}$ is a bijection to the neighbors of $v$.  Hence, the cyclic order
around $\tilde{f}^{-1}(\tilde{v})$ (which is a subset of $(\tilde{f} \circ
f)^{-1}(v)$) induces a cyclic order on the neighbors of $\tilde{v}$. This
establishes \pv.   By construction the bijection induced by $f$ between the
neighbors of $\tilde{v}$ and the neighbors of $f(\tilde{v})$ is order
preserving.

Before establishing \pe\ we must assign valid signs to the vertices of
$\widetilde{\widetilde{G}}$.  Since $\tilde{f} \circ f$ fulfills \pe, some
signs have been assigned already, and these signs are valid for $f \circ
\tilde{f}$.  We show that the same signs are valid for 
$\tilde{f}$.  Let $\tilde{\tilde{v}}, \tilde{\tilde{v}}' \in
V(\widetilde{\widetilde{G}})$ be vertices so that
$\tilde{f}(\tilde{\tilde{v}}) = \tilde{f}( \tilde{\tilde{v}}')$.  Then $f
\circ \tilde{f}(\tilde{\tilde{v}}) = f \circ \tilde{f}(\tilde{\tilde{v}}')$
and therefore $\tilde{\tilde{v}}$ and $\tilde{\tilde{v}}'$ have the same
(resp.  opposite) sign if and only if  the \ccc\ around them is the same
(resp. opposite) under $f \circ \tilde{f}$.  The order preserving bijection
induced on the neighbors of $\tilde{v}$ by $f$ shows that the sign choice is
valid for $\tilde{f}$ as well.

Let $\tilde{\tilde{e}}, \tilde{\tilde{e}}' \in E({\widetilde{\widetilde{G}}})$
be two edges that project to the same edge (say $\tilde{e}$) under
$\tilde{f}$.    By assumption $\tilde{f} \circ f$ fulfills \pe.  Since 
$\tilde{\tilde{e}}$ and $\tilde{\tilde{e}}'$ project to the same edge under
$f \circ \tilde{f}$, $\tilde{\tilde{e}}$ connects vertices of the same sign
if and only if $\tilde{\tilde{e}}'$ does. Hence
$\tilde{f}:\widetilde{\widetilde{G}} \to \widetilde{G}$ fulfills \pe.

This completes the proof of Proposition~\ref{pro:tilde-f-has-pev}.
\end{proof}

By Proposition~\ref{pro:tilde-f-has-pev} we may apply
Theorem~\ref{thm:sufficiency} to $\tilde{f}: \widetilde{\widetilde{G}} \to
\widetilde{G}$ and get an embedding of $\widetilde{G}$ into $F$, where $F
\cong S^2$ or $F \cong P^2$.

\medskip

\noindent {\bf Case One: $F \cong S^2$:}  This case corresponds to Case~(1) of
Theorem~\ref{thm:necessity}.  We need to show that $f:\widetilde{G} \to G$
fulfills \pev.  By Theorem~\ref{thm:sufficiency}~(\ref{order}) the embedding
of $\widetilde{G}$ into $S^2$ and the map $\tilde{f}$ induce the same order
around the vertices of $\widetilde{G}$.  Similarly, applying
Theorem~\ref{thm:sufficiency} to $f \circ \tilde{f}$ we obtain an embedding of
$G$ into some surface $F'$, and the embedding of $G$ into $F'$ and the map $f
\circ \tilde{f}$ induce the same order around the vertices of $G$.

To verify \pv, fix $\tilde{v} \in V(\widetilde{G})$ and denote $f(\tilde{v})$
by $v$.  By assumption, the cover $f$ is not branched.  Therefore $f$ induces
a bijection between the neighbors of $\tilde{v}$ and the neighbors of $v$.
Since the orders around $\tilde{v}$ and around $v$ are both induced from the
order around their preimages in $\widetilde{\widetilde{G}}$, this bijection is
order preserving.  Therefore, $f$ fulfills \pv.

Next we need a valid sign assignment on $V(\widetilde{G})$.  Fix $\tilde{v}
\in V(\widetilde{G})$.  By assumption, $F \cong S^2$ and therefore by
Proposition~\ref{pro:p2iff-both-signs} all the preimages of $\tilde{v}$ have
the same sign.  We assign $\tilde{v}$ that sign.  We need to verify that the
assignment is valid.  Fix $\tilde{v}, \ \tilde{v}' \in V(\widetilde{G})$ that
project to the same vertex, say $v$.  Let $D_v$ be a small disk neighborhood
of $v$ oriented so that the restriction of $f \circ \tilde{f}$ to a component
of $(f \circ \tilde{f})^{-1}(D_v)$ is orientation preserving if and only if
the sign at the corresponding preimage of $v$ is plus (this is possible since
the sign assignment on $V(\widetilde{\widetilde{G}})$ is valid for $f \circ
\tilde{f}$).  We choose an orientation for $F \cong S^2$ so that $\tilde{f}$
is orientation preserving.  It is easy to see that given $\tilde{v}$ with
$f(\tilde{v}) = v$, $f|_{D_{\tilde{v}}} :D_{\tilde{v}} \to D_v$ is orientation
preserving if and only if only if the sign at $\tilde{v}$ is plus; it follows
that the sign assignment is valid.

Finally we verify \pe.  Given $e \in E(G)$ and $\tilde{e},\ \tilde{e}'$ that
project to $e$, the signs at the endpoints of $\tilde{e}$ are the same as the
signs at the endpoints of of any edge that projects to $\tilde{e}$ under
$\tilde{f}$, and similarly for $\tilde{e}$.  \pe\ for $f$ follows from \pe\
for $\tilde{f}$ (that was established is
Proposition~\ref{pro:tilde-f-has-pev}).

\medskip

\noindent {\bf Case Two: $F \cong P^2$:}  This case corresponds to Case~(2) of
Theorem~\ref{thm:necessity}.  Let $\widehat{G}$ be the lift of $\widetilde{G}$
to $P^2$, {\it i.e.}, $\widehat{G} = \pi^{-1}(\widetilde{G})$.  (Recall the
definition of the universal cover $\pi:S^2 \to P^2$ in
Notation~\ref{notation:universal-cover}.)

Although not essential to the proof, we show that $\widehat{G}$ is connected.
Every face of $S^2$ cut open along $\widehat{G}$ is an unbranched cover of a
face of $\widetilde{G}$, which by Theorem~\ref{thm:sufficiency}~(5) (applied
to $\tilde{f}$) is a disk. Therefore the faces of $\widehat{G}$ are disks as
well and $\widehat{G}$ is connected.

By Lemma~\ref{lem:covers-of-p2-factor} $\tilde{f}$ factors through $\pi$; that
is, there exist a cover $\hat{f}:S^2 \to S^2$ so that $\tilde{f} = \pi \circ
\hat{f}$.  Then $\widehat{G} = \hat{f}(\widetilde{\widetilde{G}})$.  Therefore
the orders induced on the neighbors of every vertex of $V(\widehat{G})$ by the
embedding into $S^2$ and by $\hat{f}$ are the same, and we conclude that $f
\circ \pi$ induces that same order on the neighbors of every vertex in
$V(G)$ as $f \circ \tilde{f}$; in particular, $f \circ \pi$ induces {\it some}
order on the neighbors of every vertex in $V(G)$ and therefore fulfills \pev.

This completes the proof of Theorem~\ref{thm:necessity}.
\end{proof}

\section{Examples}
\label{sec:examples}

Our first example is very simple.  It shows a cover of a graph  with one
vertex and two edges (the bouquet of two circles). The cover is given in
Figure~\ref{fig:irregular} and is an  irregular triple cover that fulfills
\pev.  (We note that any double cover is regular.)   In that figure, we use
the following labels $a$ and $b$, where the three edges labeled $a$ project to
the same edge and the three edges labeled $b$ project to the other edge.  The
the bouquet of two circles is not shown in the figure.
\begin{example}
\label{ex:irregular}
Let $G$ be the bouquet of two circles and $\widetilde{G}$ the triple cover
given in figure~\ref{fig:irregular}; the projection $f$ is indicated by labels
and arrows.  Then $f:\widetilde{G} \to G$ is an irregular cover fulfilling
\pev.
\end{example}
\noindent Of course, the bouquet of two circles cannot be regarded as an
``interesting'' graph.  The reader can soup this example up by replacing the
edges of the bouquet of two circles by planar graphs.  However, the resulting
cover is not 3-connected.  We ask:

\begin{figure}
\psfrag{a}{$a$}
\psfrag{b}{$b$}
\centerline{\includegraphics[height=3cm]{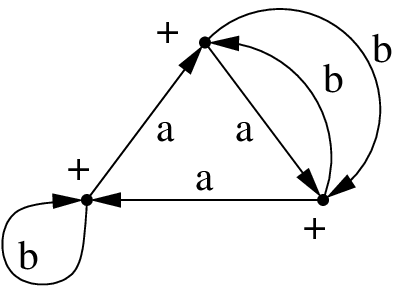} }
\label{fig:irregular}
\caption{An irregular cover can have \pev}
\end{figure}

\begin{question}
Let $f:\widetilde{G} \to G$ be a finite planar graph fulfilling \pev.  Suppose
$G$ has no 1 or 2-cycles and $\widetilde{G}$ is 3-connected, is the cover
regular?
\end{question}

All the remaining examples in this section are covers of $K_4$.  The labels
used are as follows: the vertices of $K_4$ are labeled $a$, $b$, $c$, and $d$
and vertices of the covers are labeled by the vertex they project to; edges
are not labeled.

\begin{figure}
\psfrag{a}{$a$}
\psfrag{b}{$b$}
\psfrag{c}{$c$}
\psfrag{d}{$d$}
\psfrag{+}{$+$}
\centerline{\includegraphics[height=3cm]{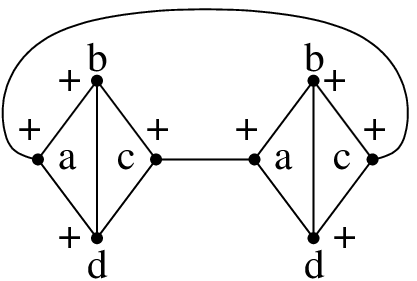} }
\caption{First of two covers of $K_4$ fulfilling \pev}
\label{fig:k4-1}
\end{figure}

\begin{figure}
\psfrag{a}{$a$}
\psfrag{b}{$b$}
\psfrag{c}{$c$}
\psfrag{d}{$d$}
\psfrag{+}{$+$}
\psfrag{-}{$-$}
\centerline{\includegraphics[height=3cm]{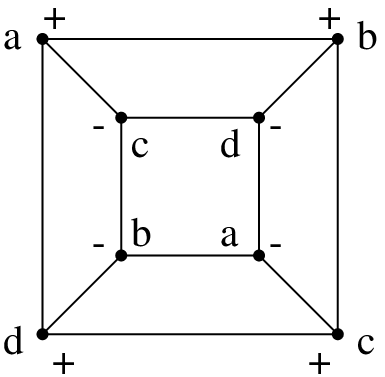} }
\caption{Second of two covers of $K_4$ fulfilling \pev}
\label{fig:k4-2}
\end{figure}

\begin{example}
\label{ex:k4}
Figures~\ref{fig:k4-1} and \ref{fig:k4-2} give two regular, unbranched, planar
double covers of  $K_4$.  The first cover yields an embedding of $K_4$ into
$S^2$ while the second cover yields an embedding of $K_4$ into $P^2$.
\end{example}
\noindent We see by inspection that both covers fulfill \pev.  It
follows from Proposition~\ref{pro:p2iff-both-signs} the first cover gives an
embedding of $K_4$ into $S^2$ while the second embeds $K_4$ in $P^2$.
Alternatively, we can see that the second cover yields an embedding into $P^2$
by observing that it has two cycles of length four that bound faces.  These
cycles project to two cycles of length two or a single cycle of length four in
$K_4$  (in fact, they project to a single cycle of length four).   By
Theorem~\ref{thm:necessity}~(5) this cycle bounds a face.  But the unique
embedding of $K_4$ into $S^2$ has only triangular faces.

\begin{figure}
\psfrag{a}{$a$}
\psfrag{b}{$b$}
\centerline{\includegraphics[height=3cm]{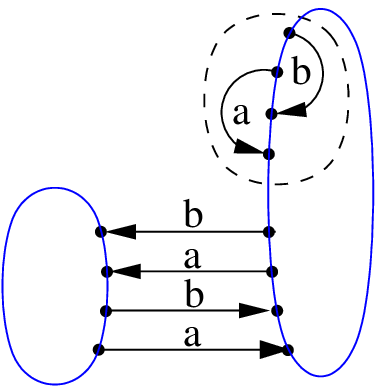} }
\caption{$\widetilde{K}_4$ in Example~\ref{ex:have-to-lift}}
\label{fig:have-to-lift1}
\end{figure}

The next example, Example~\ref{ex:have-to-lift}, is the triple cover
$f:\widetilde{K}_4 \to K_4$ given in Figure~\ref{fig:have-to-lift1}.  We
observe that with the given embedding the cover does not fulfill \pe.
This cover is not 3-connected so the embedding into $S^2$ is not
unique.  However, the only other embedding is given by ``flipping'' the disk
contained in the dashed circle.  It is easy to see that it too does not
fulfill \pe.  

We seek a planar cover $\widetilde{\widetilde{K}}_4
\to \widetilde{K}_4$ so that the composition 
$\widetilde{\widetilde{K}}_4 \to K_4$ fulfills \pev.
Although covers of this type may have arbitrarily high degree, by
Theorem~\ref{thm:necessity} if such cover exists, then a double cover with
this property exists
as well.  Moreover, Theorem~\ref{thm:necessity} tells us how to find this
cover: let $\pi:S^2 \to P^2$ be the universal cover; embed $\widetilde{K}_4$
in  $P^2$ and lift to $S^2$; this is the cover we need.

\begin{figure}
\psfrag{a}{$a$}
\psfrag{b}{$b$}
\centerline{\includegraphics[height=3cm]{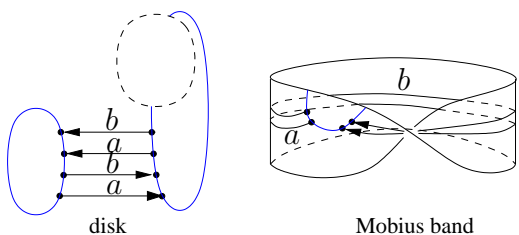} }
\caption{Embedding $\widetilde{K}_4$ in $P^2$}
\label{fig:have-to-lift2}
\end{figure}

\begin{figure}
\psfrag{a}{$a$}
\psfrag{b}{$b$}
\centerline{\includegraphics[height=3cm]{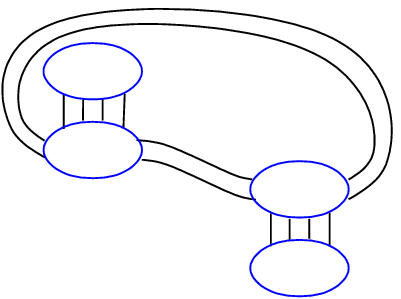} }
\caption{The lift to the universal cover $S^2$}
\label{fig:have-to-lift4}
\end{figure}

\begin{example}
\label{ex:have-to-lift}
The cover given in Figure~\ref{fig:have-to-lift1} does not fulfill \pev\ for
any embedding of $\widetilde{G}$ into $S^2$.  However, the embedding of $\widetilde{K}_4$
into $P^2$ given in  Figure~\ref{fig:have-to-lift2} (where we view $P^2$ as a
disk union a M\"obius band) has a double cover (given in
Figure~\ref{fig:have-to-lift4}) that embeds in $S^2$ and fulfills \pev.
\end{example}

\begin{figure}
\psfrag{a}{$a$}
\psfrag{b}{$b$}
\centerline{\includegraphics[height=3cm]{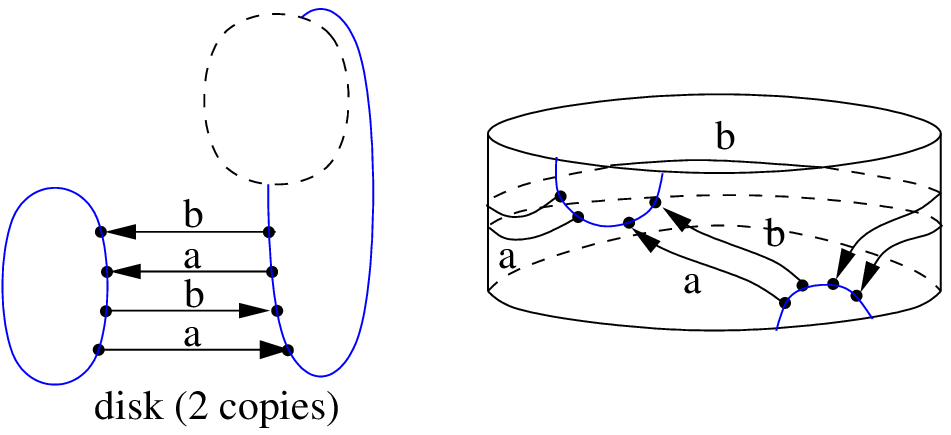} }
\caption{The lift to the universal cover $S^2$}
\label{fig:have-to-lift3}
\end{figure}

Figure~\ref{fig:have-to-lift4} was constructed by taking two copies of the
disk in Figure~\ref{fig:have-to-lift2} and an annulus that double covers the
M\"obius band in  Figure~\ref{fig:have-to-lift2}; see
Figure~\ref{fig:have-to-lift3}.  Pasting the disks to the annulus in
Figure~\ref{fig:have-to-lift3} gives Figure~\ref{fig:have-to-lift4}.

Our final example is far more disturbing in nature.  We give a cover
$f:\widetilde{K}_4 \to K_4$ that does not fulfill \pev\ (quite similar to
Example~\ref{ex:have-to-lift}) and {\it conjecture} that for any finite
planar cover $\widetilde{\widetilde{K}}_4 \to \widetilde{K}_4$, the
composition $\widetilde{\widetilde{K}}_4 \to K_4$ does not fulfilling \pev.

\begin{figure}
\centerline{\includegraphics[height=3cm]{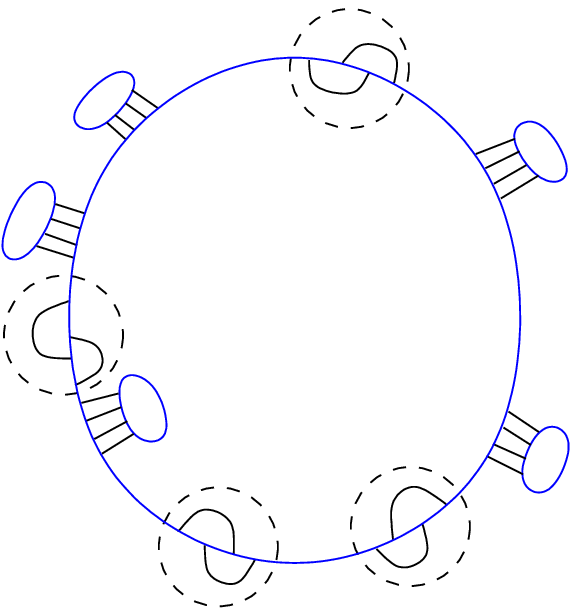} }
\caption{Can we obtain \pev\ by passing to a double cover?}
\label{fig:no-properties}
\end{figure}

\begin{cnj}
\label{ex:no properties}
Consider the cover $\widetilde{K}_4$ given in Figure~\ref{fig:no-properties}.
Let $\widetilde{\widetilde{K}}_4 \to \widetilde{K}_4$ be a finite planar
cover.  Then the composition $\widetilde{\widetilde{K}}_4$ does not fulfill
\pev.
\end{cnj}

\begin{figure}
\centerline{\includegraphics[height=1.5cm]{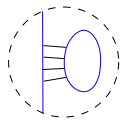} }
\caption{Fixing the disks in Figure~\ref{fig:no-properties}}
\label{fig:fixme}
\end{figure}

\begin{figure}
\centerline{\includegraphics[height=3cm]{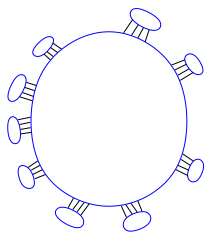} }
\caption{The fixed cover}
\label{fig:fixed}
\end{figure}

We end this paper by describing two ways of dealing with the cover
$\widetilde{K}_4$ given in Figure~\ref{fig:no-properties}.  First, note that
if the four disks enclosed in dashed circles in Figure~\ref{fig:no-properties}
are replaced by the disks given in Figure~\ref{fig:fixme} (and some of the
``mushrooms'' are reflected), we obtain the cover shown in
Figure~\ref{fig:fixed} that does fulfill \pev. The work in this paper suggests
that topological ``cut-and-paste'' techniques such as this could be very
useful.

Next, recall Figures~\ref{fig:have-to-lift1} and \ref{fig:have-to-lift2} where
we replaced a disk with a M\"obius band.  Similarly, we can replace the four
disks in Figure~\ref{fig:no-properties} with four M\"obius bands, obtaining an
embedding of $\widetilde{K}_4$ into a non-orientable surface $S$; since
replacing a disk with a M\"obius band lowers the Euler characteristic by 1,
$\chi(S) = -2$.  Just like ${P}^2$, $S$ has an oriented double cover (say
$\widetilde{S}$).  The Euler characteristic is multiplicative under unbranched
covers, and we conclude that $\widetilde{S}$ is the surface of genus 3.  (In
general, if we replace $n > 0$ disks with M\"obius bands the oriented double
cover will have genus $n-1$.)  The reader can verify that the proof of
Theorem~\ref{thm:sufficiency} is valid, and gives a surface $F$, an embedding
$G \subset F$,  and a cover $f':\widetilde{S} \to F$ that extends $f$.  As in
Lemma~\ref{pro:F-covered-bySorP}, $\chi(F) \geq -2$.  (In general, $\chi(F)
\geq 2-n$.)  This does not give an embedding of $G$ into $P^2$, but does give
an embedding of $G$ into a surface with some control over its Euler
characteristic.

\end{document}